\documentclass[12pt]{amsart}
\title{On a generalization of the plank problem}
\author[Zs. P\'ales]{By Zs. P\'ales}
\address{Institute of Mathematics, L. Kossuth University, H-4010 Debrecen, Pf. 12, Hungary.}
\email{pales@math.klte.hu}
\date{}

\begin{document}
%\large
\begin{flushright}
{\sl General Inequalities 6}, (Oberwolfach, 1990), \\
(ed.\ W.\ Walter), Birkh\"{a}user Verlag, \\
Basel-Boston-Stuttgart, 1992, pp.\ 473-476.\\[3mm]
\end{flushright}

\maketitle

\def\R{{\mathbb R}}
\def\Int{\int\limits}
\let\al=\alpha
\let\ep=\varepsilon
\newtheorem{thm}{Theorem}

A {\it strip} or a {\it plank} $S$ in $\R^n$ is a closed set bounded by
two parallel hyperplanes.  The distance of these hyperplanes is called
the {\it width} of $S$.  The {\it minimal width} of a convex closed set
$K$ is the minimal width of a strip containing $K$.

The following theorem was conjectured by A. Tarski in 1932 and proved
by T. Bang [2] in 1951:

{\it If a closed convex set $K$ in $\R^n$ is covered by a finite number
of strips, then the sum of their widths is greater than or equal to the
minimal width of $K$.}

This result has recently been generalized to Banach spaces by K. Ball
[1].

If $n=2$ and $K$ is the unit disc then there is an extremally simple
proof for the above result.

Assume that the unit disc in $\R^2$ is covered by strips $S_1,\dots,S_k$
with widths $d_1,\dots,d_k$. Without loss of generality we can also
assume that both bounding lines of the strips intersect the unit circle.
Now we consider the unit sphere in $\R^3$ and to each strip $S_i$ in
$\R^2$ we construct a three dimensional strip $S_i^*$ which is of width
$d_i$ and intersects the $xy$-plane in $S_i$. Since $S_1,\dots,S_k$
cover the unit disc, hence $S_1^*,\dots,S_k^*$ cover the unit sphere.
The area of the inersection of the unit sphere and the strip $S_i^*$ is
$2\pi d_i$ independently of the position of of the $i$-th strip. (This
is a well known fact from calculus, already discovered by Archimedes.)
Thus the sum of these areas exceeds the area of the unit sphere, i.e.
$$
  \sum_{i=1}^k 2\pi d_i \ge 4\pi \qquad \Longrightarrow \qquad
  \sum_{i=1}^k d_i \ge 2,
$$
which was to be proved.

We can interpret this proof in the following way: If $S$ is a subset of
the disc then we project it up to the sphere, measure the area of the
projection and call this number the $\mu$ measure of $S$.  Then the
$\mu$ measure of a the intersection of a strip and the disc is the width
of the strip times $2\pi$.  Then the statement is a simple consequence
of the subadditivity of $\mu$.  In what follows, we generalize this idea
and extend the result discussed above.

An {\it angular domain} in $\R^2$ is a closed convex set $D$ bounded by
two halflines. The {\it angle} of $D$ is the angle closed by the
bounding halflines. The {\it vertex} of $D$ is the common endpoint of
these two halflines.

\begin{thm}
Let two concentric circles $k$ and $K$ be given on
the plane with radii $r$ and $R$, $r<R$.  Assume that the disc
bounded by $k$ is covered by angular domains whose vertices are within $K$.
Then the sum of the angles of these angular domains is greater than or equal to
the view angle of $k$ from an arbitrary point of $K$.
\end{thm}

{\bf Remark.} This result was proposed as a problem by the author on the
1985 M. Schweitzer competition (see [3]).

{\it Proof.} Denote by $O$ the common center of the circles and by
$D_1,\dots,D_k$ the given angular domains with angles $\al_1,\dots,\al_k$.
An angular domain $D$ will be called {\it regular} if the vertex of $D$ is
on $K$ and both bounding halflines of $D$ intersect $k$.  Without
loss of generality, we can assume that $D_1,\dots,D_k$ are regular
domains.

The idea of the proof is the following: We construct a rotation
invariant nonnegative measure $\mu$ on the closed disc $T$ bounded by
$k$ such that the measure of the intersection of $D$ and $T$ is $\al$,
where $D$ is an arbitrary regular angular domain with angle $\al$. Having
such a measure we can give a one line proof for the theorem:
$$
  \sum_{i=1}^k\al_i=\sum_{i=1}^k \mu(D_i\cap T)\ge
  \mu\biggl(\bigcup_{i=1}^k (D_i\cap T)\biggr)=\mu(T),
$$
and observe that $\mu(T)$ is exeactly the view angle of $k$ from any
point of $K$.

Now we construct the desired measure. Let $P\in T$ be an arbitrary
point, and denote by $\rho$ the distance of $P$ and $O$. Then define
$$
  F(P)=f(\rho)={1\over\pi}\cdot{1\over R^2-\rho^2}\cdot
       \sqrt{R^2-r^2\over r^2-\rho^2}.
$$
If $S$ is a Lebesgue measurable subset of $T$ then let
$$
  \mu(S)=\Int_S F(P)dP.
$$
Obviously $\mu$ is a rotation invariant nonnegative measure on $T$.
To prove the key property of $\mu$, let $D$ be a regular angular domain
with vertex $A$ and angle $\al$. Then we want to show $\mu(D)=\al$.
Without loss of generality we can assume that one of the bounding
halflines of $D$ is tangent to $k$ at the point $Q$. (In the general
case $D$ can be obtained as  the difference of two such angular domains.)
Denote by $\ep$ the angle $OAQ<$ and by $d$ the (signed) distance of the
other bounding halfline from $O$. (This distance is positive if $O$ is
outside $D$, and negative if $O$ is inside $D$.) Then
$d=R\sin(\ep-\al)$. Using succesive integration, we obtain
\begin{eqnarray*}
  \mu(D)  &  = &\Int_D F(P)dP
             =\Int_d^r\Int_{-\arccos(d/\rho)}^{\arccos(d/\rho)}
                     f(\rho)\rho d\varphi d\rho                    \\
          &  = &\Int_d^r 2f(\rho)\rho \arccos(d/\rho)d\rho         \\
          &  = &\Int_{R\sin(\ep-\al)}^r 2f(\rho)\rho
                  \arccos\left({R\sin(\ep-\al)\over\rho}\right)d\rho.
\end{eqnarray*}
Thus we have to show that
$$
    \Int_{R\sin(\ep-\al)}^r 2f(\rho)\rho
         \arccos\left({R\sin(\ep-\al)\over\rho}\right)d\rho=\al
$$
for all $0 \le \al \le 2\ep$. Substituting the new variable
$t=R\sin(\ep-\al)$ this reduces to
$$
  \Int_t^r 2f(\rho)\rho \arccos(t/\rho)d\rho=\ep-\arcsin(t/R),
$$
for $-r\le t\le r$. This latter equation is obviously valid for $t=r$,
thus it suffices to show that the derivatives of both sides with respect
to $t$ are identical, i.e.
$$
  \Int_t^r{2f(\rho)\rho\over\sqrt{\rho^2-t^2}}d\rho=
               {1\over\sqrt{R^2-t^2}},\qquad -r<t<r.
$$
However
\begin{eqnarray*}
  \Int_t^r{2f(\rho)\rho\over\sqrt{\rho^2-t^2}}d\rho
      & = &\Int_t^r{2\over\pi}\cdot{\rho\over R^2-\rho^2}\cdot
            \sqrt{R^2-r^2\over(r^2-\rho^2)(\rho^2-t^2)}d\rho     \\
      & = &\left[{2\over\pi}\cdot{1\over\sqrt{R^2-t^2}}\cdot
          \arctan\sqrt{{R^2-r^2\over R^2-t^2}\cdot
            {\rho^2-t^2\over r^2-\rho^2}}\right]_{\rho=t}^{\rho=r} \\
      & = &{1\over\sqrt{R^2-t^2}}.
\end{eqnarray*}
Thus the proof is complete.

{\bf Remark.} When $n=2$ and $K$ is the unit disc, then the statement of
the plank problem can easily be derived from our theorem.  Denote the
unit disc by $T$ and assume that it is covered by strips $S_1,\dots,S_k$
(whose bounding lines intersect $T$).  Take a concentric circle $K$ with
radius $R$, where $R$ is sufficiently large.  Assume that the
two bounding  lines of $S_i$ intersect $K$ in $A_i,B_i$ and in
$C_i,D_i$. We choose the notation such that $S_i$ is
covered by the two regular angular domains $A_iB_iD_i<$ and $B_iD_iC_i<$.
Denote by $\al_i'$ and $\al_i''$ their angle and by $d_i$ the width of
$S_i$. Then we have
$$
  {d_i\over 2R-2}\ge \tan\al_i'\ge\al_i',\qquad {d_i\over 2R-2}\ge
       \tan\al_i''\ge\al_i''.
$$
Thus the theorem yields
$$
    2\sum_{i=1}^k{d_i\over 2R-2}\ge\sum_{i=1}^k(\al_i'+\al_i'')\ge
        2\ep \ge 2\sin\ep\ge{2\over R}.
$$
Now taking the limit $R\to\infty$ we obtain the statement.

\end{document}